\documentclass[12pt]{amsart}

\usepackage{graphicx}
\usepackage{psfrag}

\newtheorem{thm}{Theorem}

\newtheorem{defi}{Definition}

\begin{document}

\title{Prime Decomposition of Knots in Lorenz-like Templates}

\author{Mike Sullivan}

\address{University of Texas at Austin}

\date{\today}

\begin{abstract}
In \cite{W} R. F. Williams showed that all knots in the Lorenz template are 
prime.
His proof included the cases where any number of positive twists were added
to one of the template's branches. However \cite{W} does give an example of
a composite knot in a template with a single negative twist. 
Below we will show that in all the negative 
cases composite knots do exist, and give a mechanism for producing many 
examples. This problem was cited in a list of problems in dynamics in 
\cite[problem 4.2]{BB}.
\end{abstract}

\maketitle

\section{Introduction}

\subsection{Flows, Knots and Templates}

The periodic orbits of a flow in a 3-manifold form knots. These knots and
how they are linked have been studied with the aid of {\em templates} or
{\em knots holders}, i.e. 2-dimensional branched manifolds with semi-flows 
\cite{BW}.

\begin{figure}[htb]
	\begin{center} 
        \includegraphics[height=2.5in,width=2.5in]{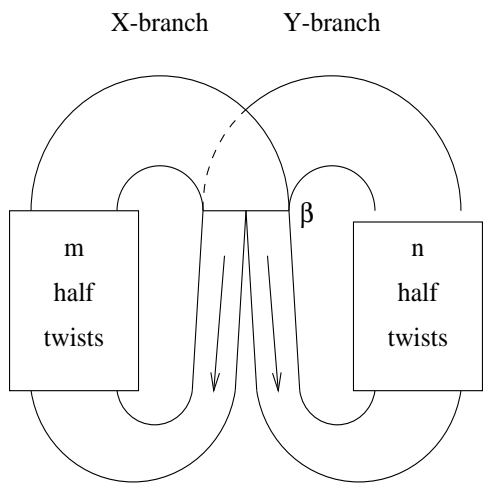}
	\end{center}
	\caption{Lorenz-like template,$L(m,n)$, defined}
	\label{fig1.eps}
\end{figure}

\begin{defi}
	The $L(m,n)$ templates are shown in Figure~\ref{fig1.eps}.
	The semi-flow comes down from the branch line, $\beta$, forming
	two branches. The X-branch has $m$ half twists, while the
	Y-branch has $n$. The sign convention for the twists is that 
	left-handed twists are positive and right-handed ones are negative.
	The two branches meet tangentially at the branch line.
	
	$\tilde{L}(m,n)$ is the mirror image of $L(m,n)$.

	The symbols defined above will also be used to represent the set of 
        knots formed by the periodic orbits in the semi-flow of the template. 

	$L(0,0)$ is the {\em Lorenz template \cite{bw}}. The remainder 
        will be called {\em Lorenz-like templates}. We also note that 
        $L(m,n) = L(n,m)$.
\end{defi}

% @@@@@@@@@@@@@@@@@@@@@@@ FIGURE @@@@@@@@@@@@@@@@@@@@@@@@@@@@@@@@@@@@@@@@@@
% @@@@@@@@@@@@@@@@@@@@@@@@@@@@@@@@@@@@@@@@@@@@@@@@@@@@@@@@@@@@@@@@@@@@@@@@@@

% @@@@@@@@@@@@@@@@@@@@@@@ FIGURE @@@@@@@@@@@@@@@@@@@@@@@@@@@@@@@@@@@@@@@@@@
\begin{figure}[h!]
	\begin{center} 
        \includegraphics[height=2in,width=5in]{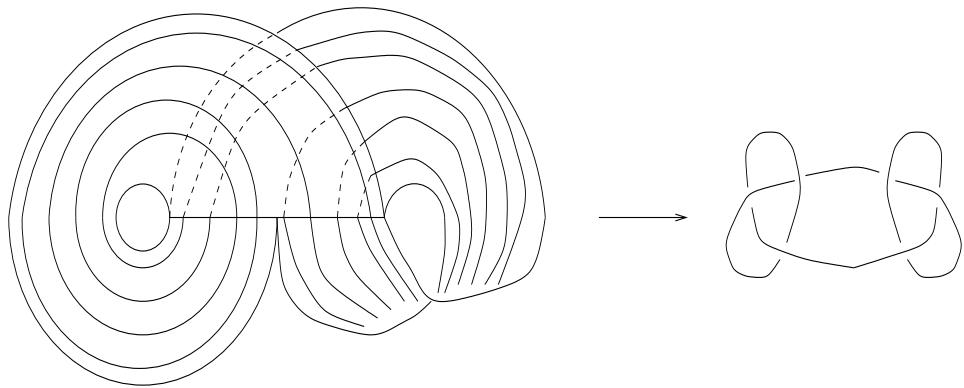}
	\end{center}
	\caption{A square knot in $L(0,-1)$}
	\label{fig2.eps}
\end{figure}
% @@@@@@@@@@@@@@@@@@@@@@@@@@@@@@@@@@@@@@@@@@@@@@@@@@@@@@@@@@@@@@@@@@@@@@@@@@

\subsection{Prime and Composite Knots}

\begin{defi}
A knot $k \subset S^{3}$ is {\em composite} if there exists a tame
sphere $S^{2}$ such that $S^{2} \cap k$ is just two points, $p$ and $q$, and 
if $\gamma$ is any arc on $S^{2}$ joining $p$ to $q$, then the knots
	\begin{eqnarray*}
		& k_{1} = \gamma \cup (k \cap \mbox{outside of } S^{2}) & \mbox{and}\\
		& k_{2} = \gamma \cup (k \cap \mbox{inside of } S^{2}), &
	\end{eqnarray*}
are nontrivial,(i.e. not the unknot). We call $k_{1}$ and $k_{2}$ 
{\em factors} of $k$ and write
	\begin{eqnarray*}
		k = k_{1} \# k_{2}.
	\end{eqnarray*}
If a knot isn't composite then it is {\em prime}.
\end{defi}

\begin{figure}[htb]
	\begin{center} 
        \includegraphics[height=2in,width=5in]{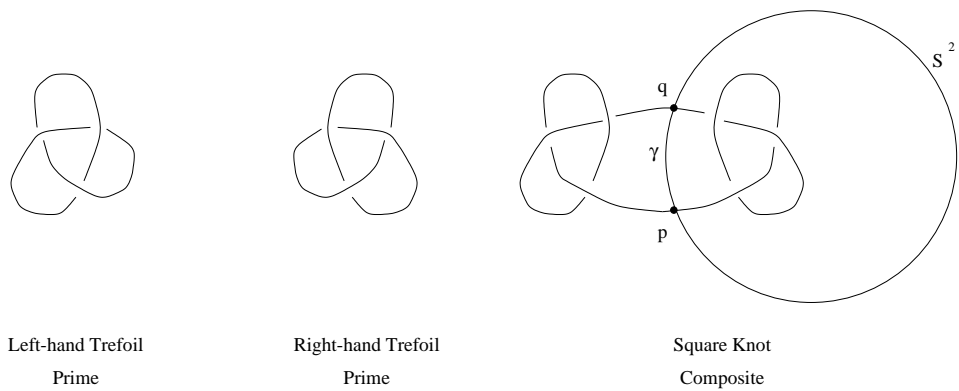}
	\end{center}
	\caption{A connected sum}
	\label{fig3.eps}
\end{figure}

Figure~\ref{fig3.eps} gives an example.

% @@@@@@@@@@@@@@@@@@@@@@@ FIGURE @@@@@@@@@@@@@@@@@@@@@@@@@@@@@@@@@@@@@@@@@@
% @@@@@@@@@@@@@@@@@@@@@@@@@@@@@@@@@@@@@@@@@@@@@@@@@@@@@@@@@@@@@@@@@@@@@@@@@@

Finally we state the all important theorem due to Schubert \cite{S,K}: 
\begin{center}
{\em Knots can be factored uniquely into primes, up to order.}
\end{center}

\section{Knots in $L(m,n)$}

\begin{thm} \label{thm_Ln-2inLn}
	As sets of knots, $L(0,n-2) \supseteq L(0,n)$, for all integers $n$.
\end{thm}

This may be restated as 

\begin{center}
 $... \supseteq L(0,-2) \supset L(0,0) \supseteq L(0,2) \supseteq ...$,

\vspace{.14in}

 $... \supseteq L(0,-3) \supseteq L(0,-1) \supset L(0,1) \supseteq L(0,3) 
\supseteq ...$,
\end{center}

\noindent where the orientable and nonorientable cases are listed separately.
We have also noted that two of the set inclusions are proper. Lorenz knots
are positive braids \cite{bw}, but $L(0,-2)$ has a negative braid. 
(To construct one, trace out the periodic orbit that wraps around the X-branch
once and the Y-branch three times.)
Knots in $L(0,1)$ are prime while  $L(0,-1)$ has composite knots. 
It is not known which, if any, of the other set inclusions are proper.

\begin{figure}[htb]
	\begin{center}
        \includegraphics[height=3.5in,width=5in]{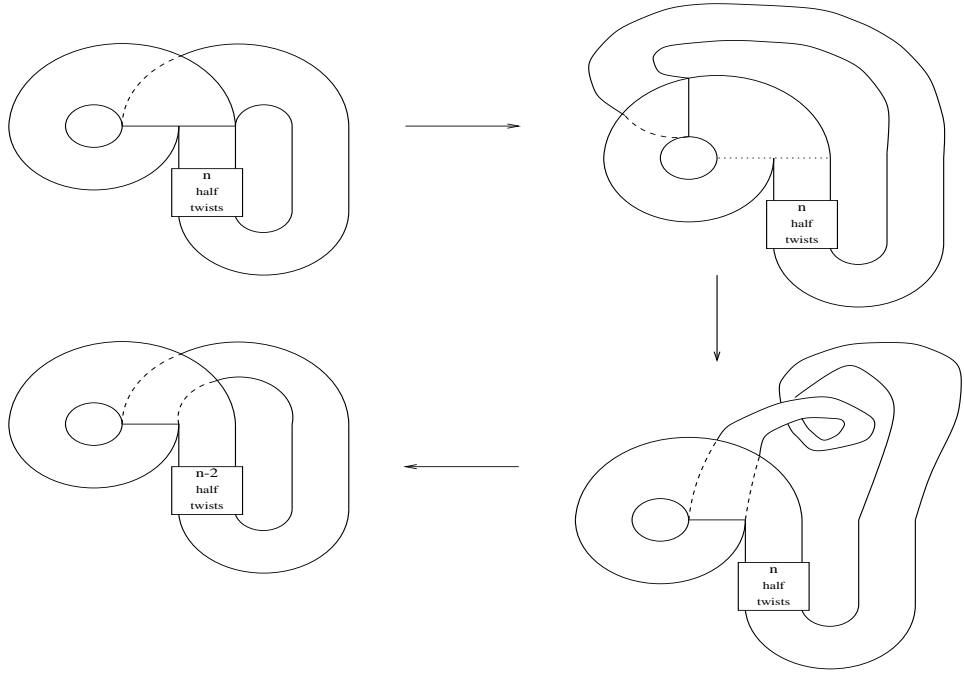} 
	\end{center}
	\caption{Proof of Theorem~\ref{thm_Ln-2inLn}}
	\label{fig4.eps}
\end{figure}

\begin{proof}
Recall the equation ``writhe = twist''. We apply this in 
Figure~\ref{fig4.eps}. The illustration shows how to take a knot embedded 
in $L(0,n)$ over to a knot of the same knot type in $L(0,n-2)$.
\end{proof}

% @@@@@@@@@@@@@@@@@@@@@@@ FIGURE @@@@@@@@@@@@@@@@@@@@@@@@@@@@@@@@@@@@@@@@@@
% @@@@@@@@@@@@@@@@@@@@@@@@@@@@@@@@@@@@@@@@@@@@@@@@@@@@@@@@@@@@@@@@@@@@@@@@@@

\begin{thm} \label{thm_main}
Any knot formed by the connected sum of a knot in $L(0,2)$ and
a knot in $\tilde{L}(0,2)$ is in $L(0,-2)$.
\end{thm}

\begin{proof}
The proof again is pictorial. We start with $L(-2,0)$, which is 
identical to $L(0,-2)$. By cutting along certain orbits and gradually 
deforming the template we will make the knots in question become clearly
visible.

Most of the steps in Figure~\ref{a-w}(a-w) should be self evident. 
Figures (a) to (n) show that the subtemplate of (n) lives in $L(-2,0)$. 
In going from figure (j) to (k) we have thrown away part of the template 
to help make the remaining steps clearer. The template of (n) is shown in (o) 
with lines to cut along. This divides (o) into two parts, (p) and (t). To see 
that any pair of orbits on (p) and (t) respectively can be summed on (n) just 
study the example in Figure~\ref{fig6.eps}. Following from (p) to (s) we get 
$\tilde{L}(0,2)$. Likewise, (t) to (w) gives $L(0,-2)$.

In Figure~\ref{fig6.eps} we show a composite knot on (n) which is the square 
knot. In Figure~\ref{bigknot.eps} we display this knot on the original form 
of the template. The reader can check this by reversing all the steps while 
carrying the knot back. 
\end{proof}

\pagebreak

a)      \hspace{2.5in}          b)

\includegraphics[height=1in,width=2in]{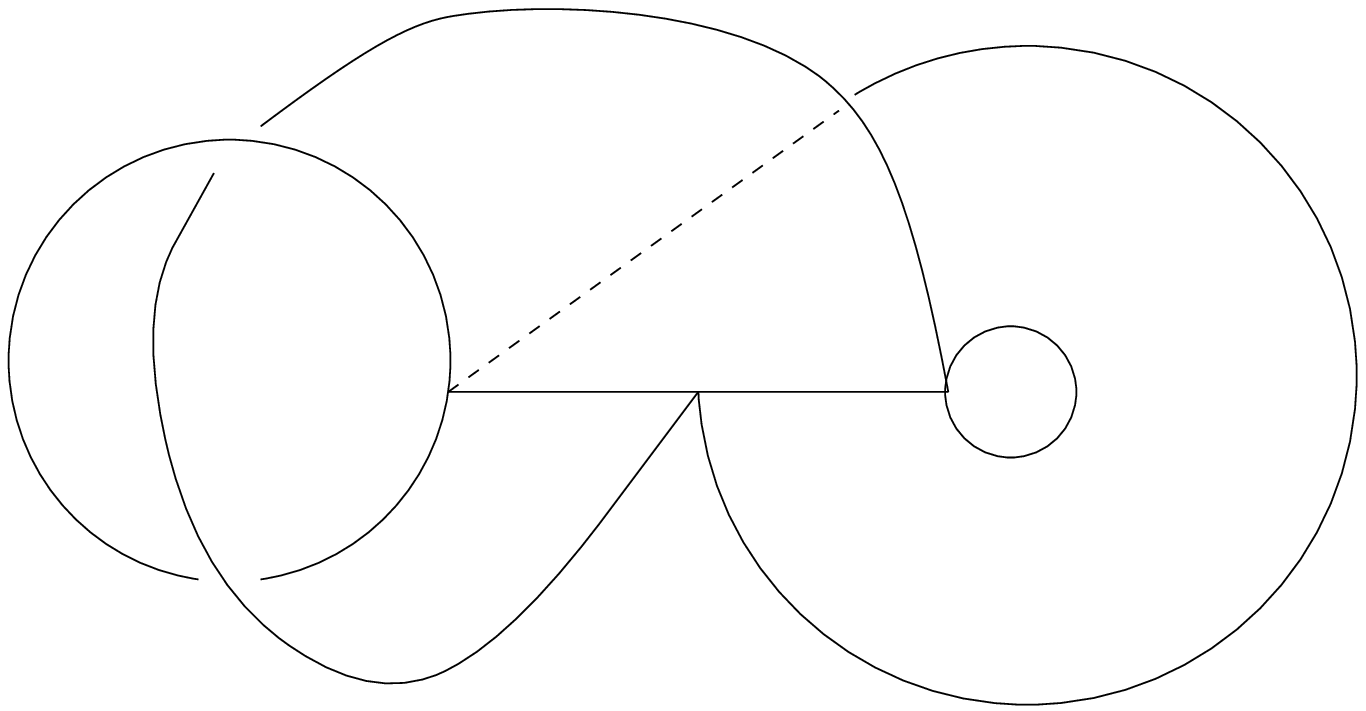} \hspace{.25in}
\includegraphics[height=1in,width=2in]{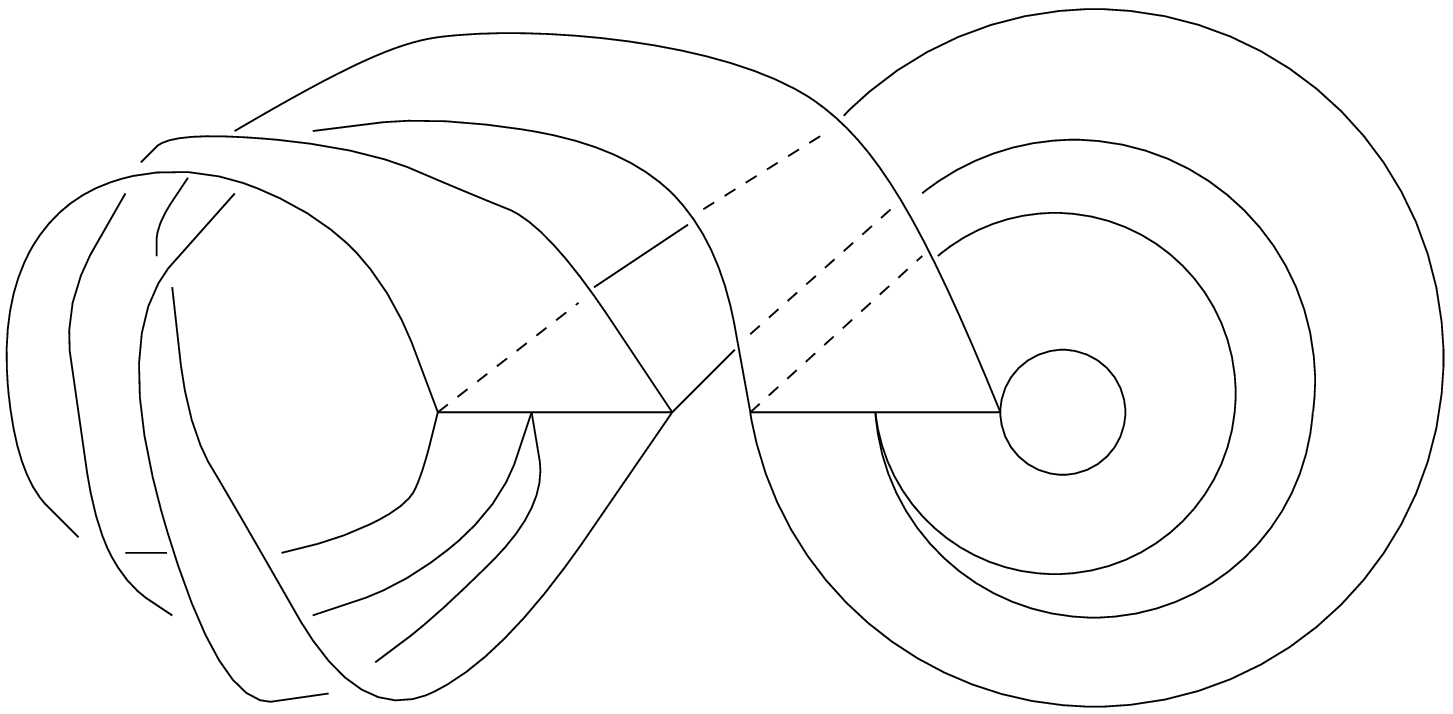}

c)      \hspace{2.5in}          d)

\includegraphics[height=1.5in,width=2in]{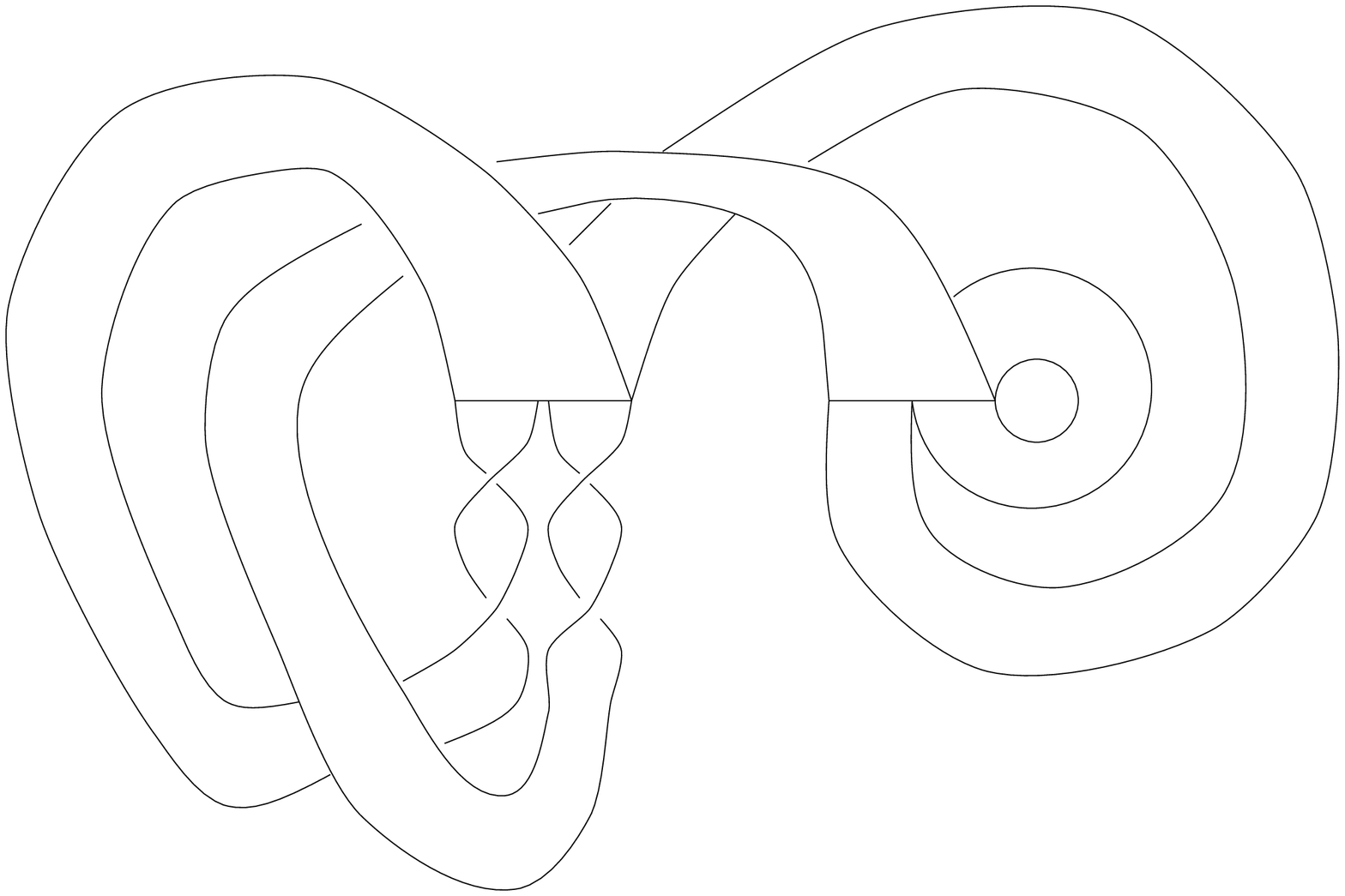} \hspace{.25in}
\includegraphics[height=1.5in,width=2in]{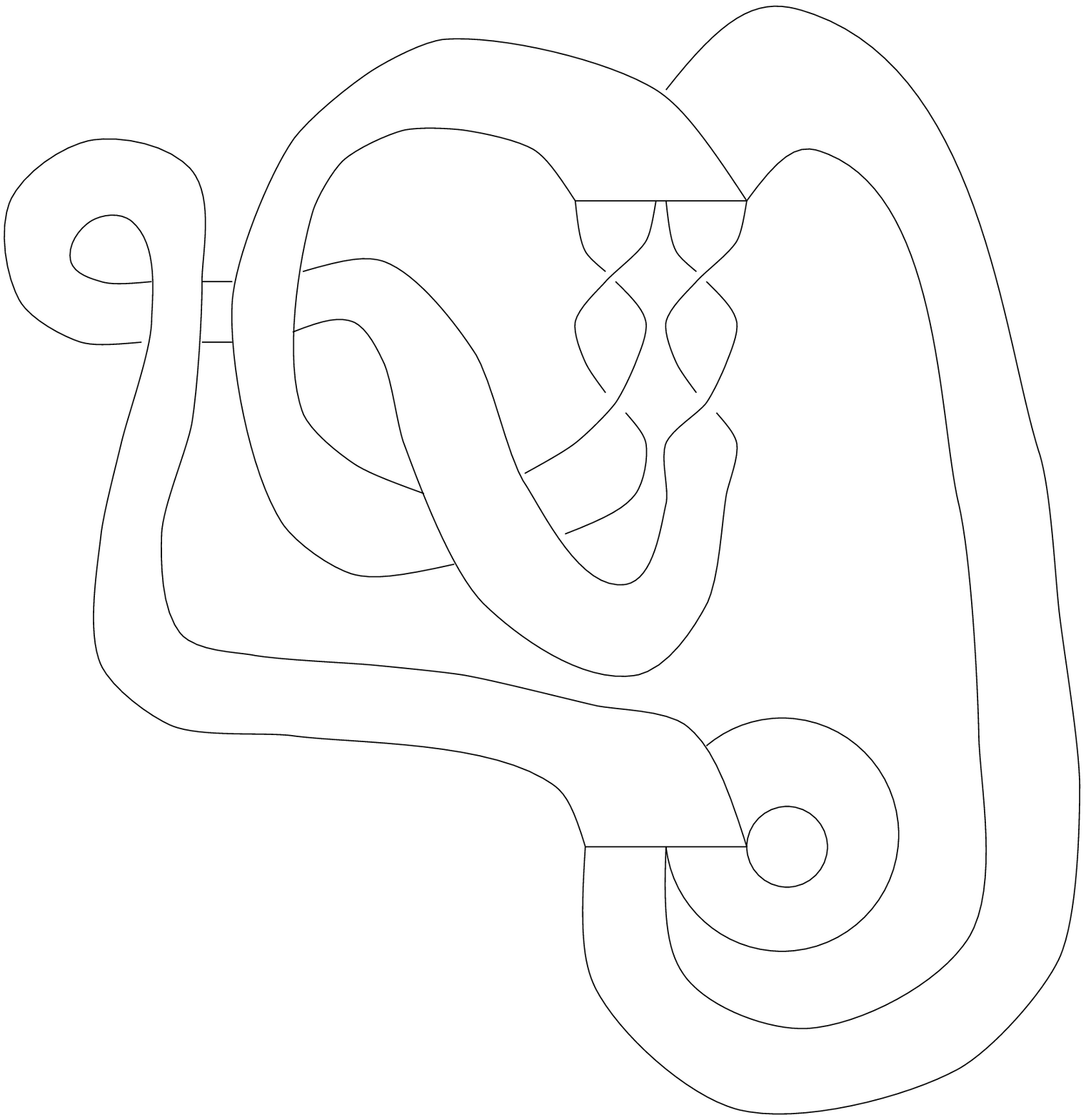} 

e)      \hspace{2.5in}          f)

\includegraphics[height=1.5in,width=2in]{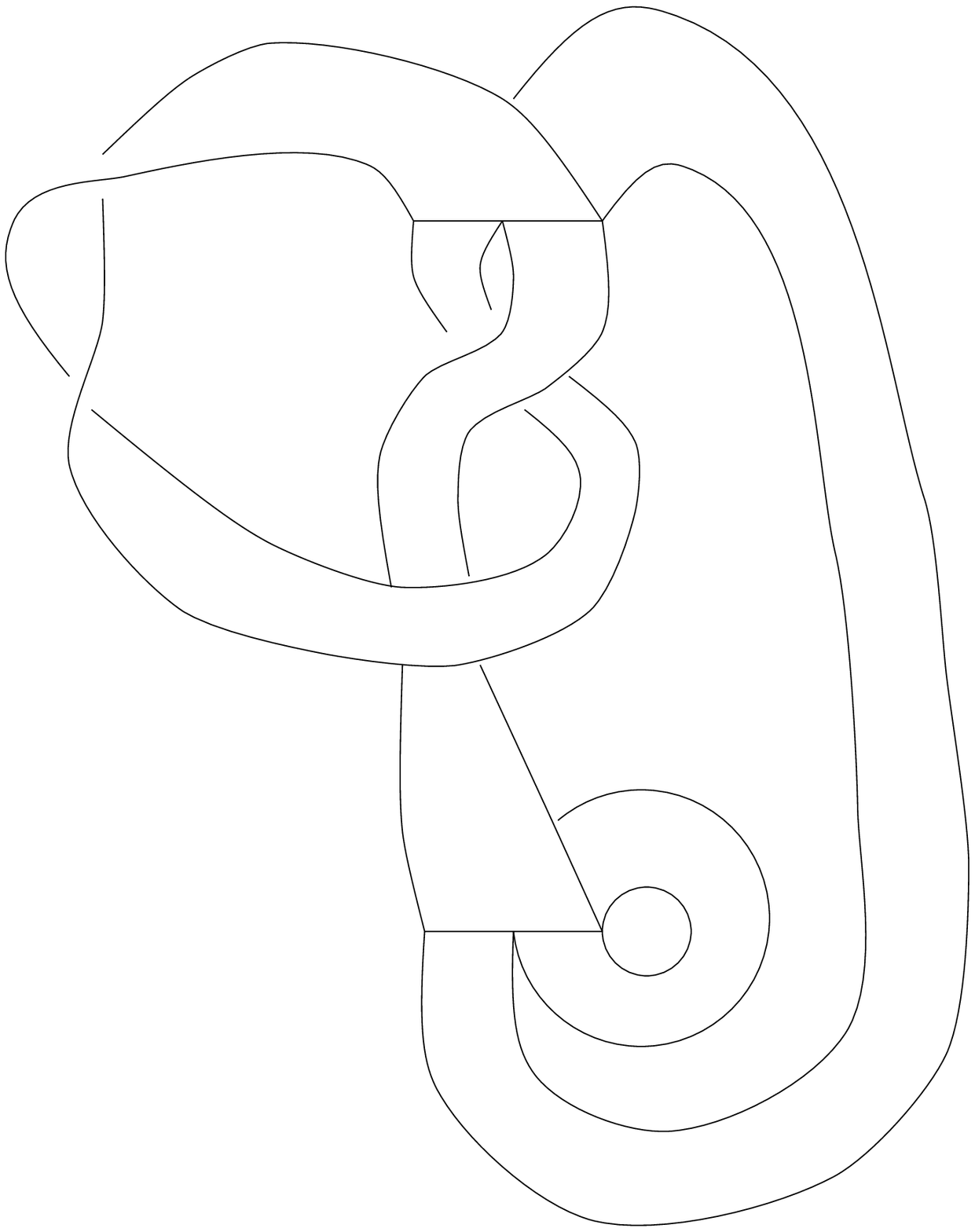} \hspace{.25in}
\includegraphics[height=1.5in,width=2in]{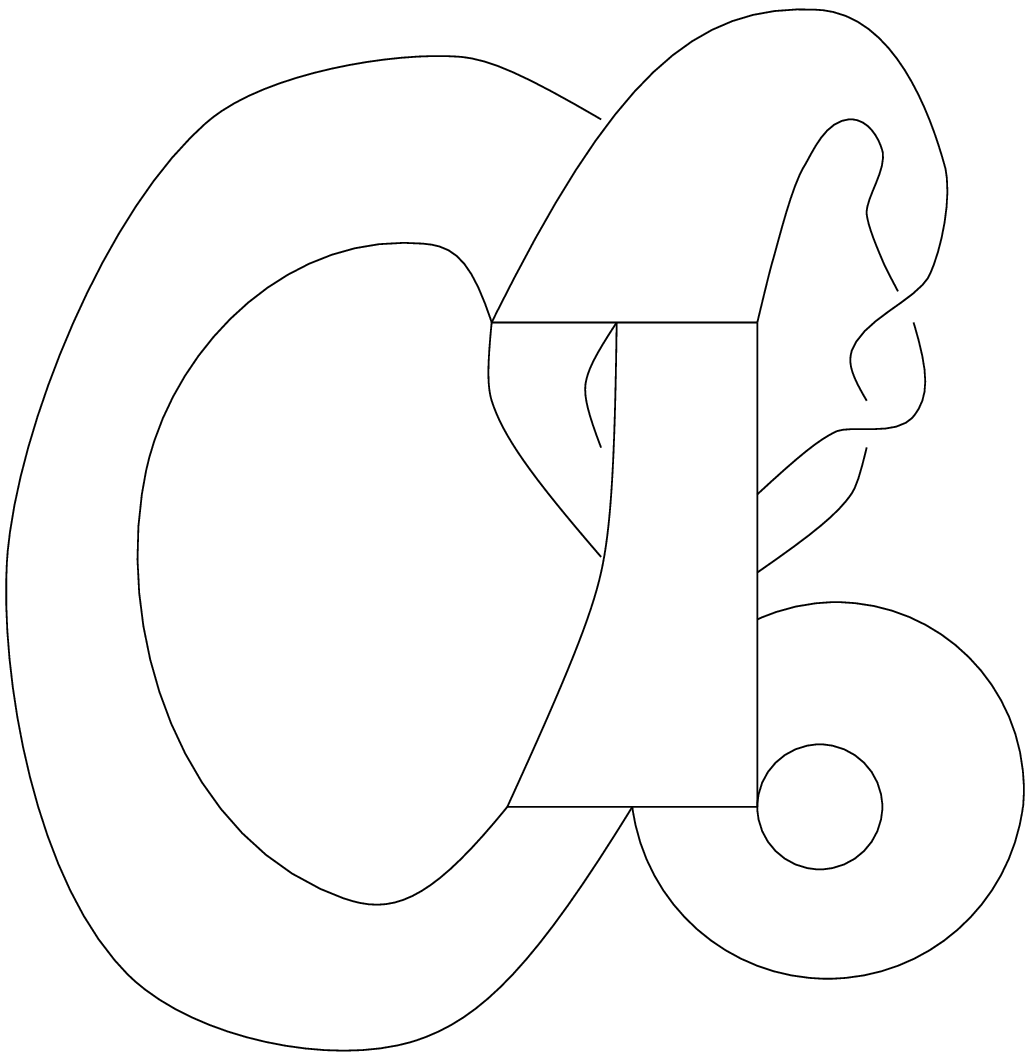} 

g)      \hspace{2.5in}          h)

\includegraphics[height=1.5in,width=2in]{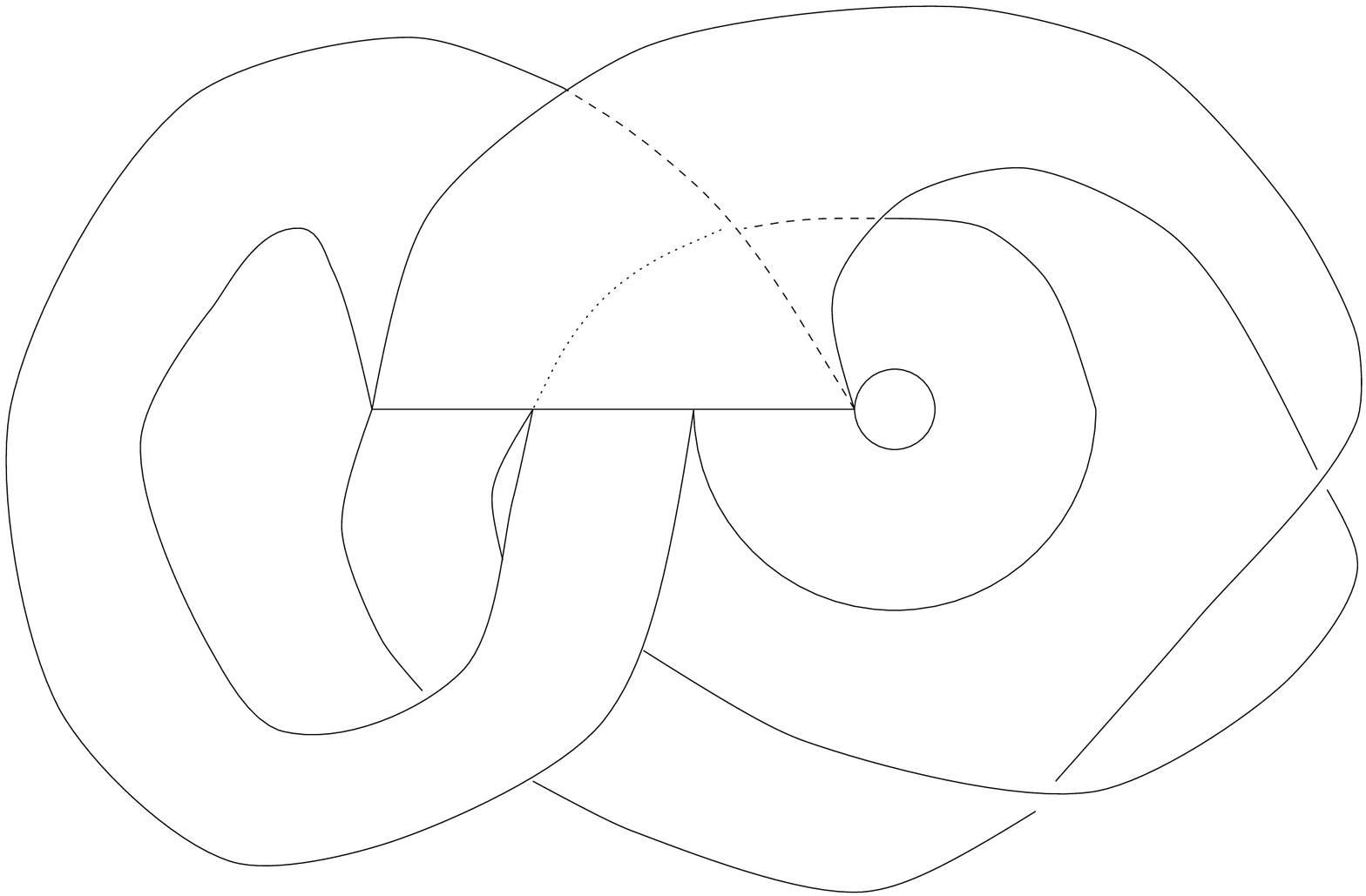} \hspace{.25in}
\includegraphics[height=1.5in,width=2in]{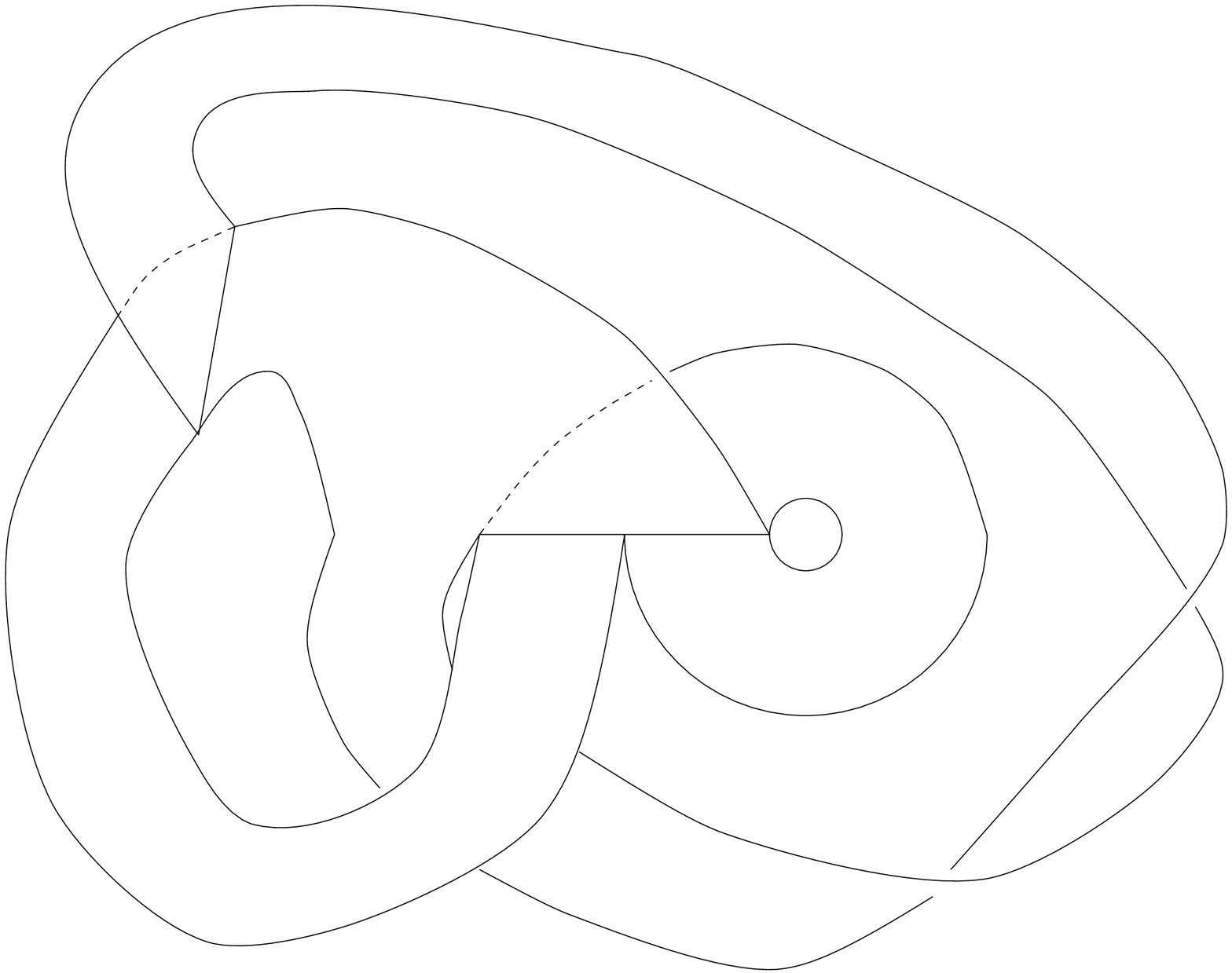}

i)      \hspace{2.5in}          j)

\includegraphics[height=1.0in,width=2in]{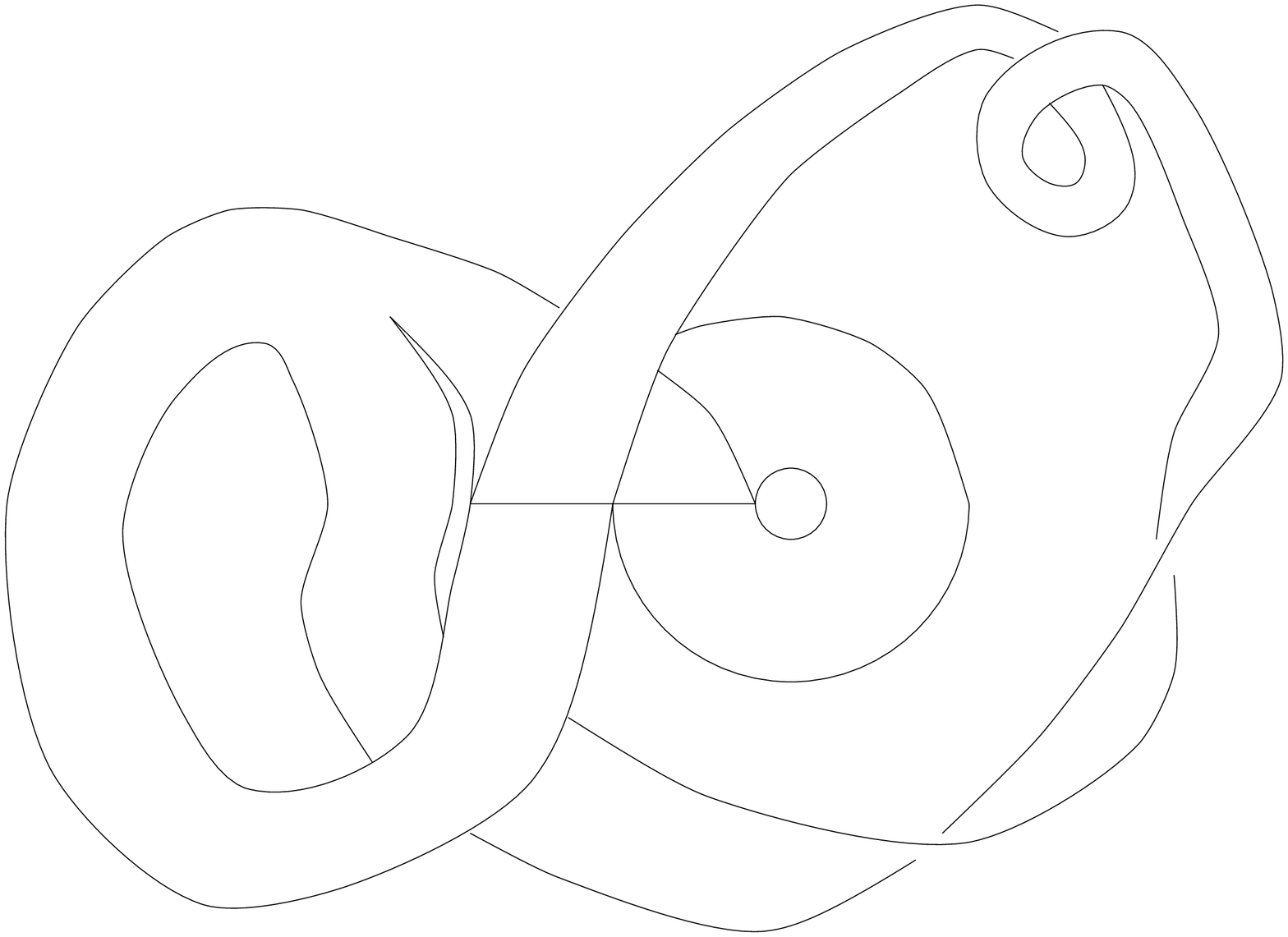} \hspace{.25in}
\includegraphics[height=1.0in,width=2in]{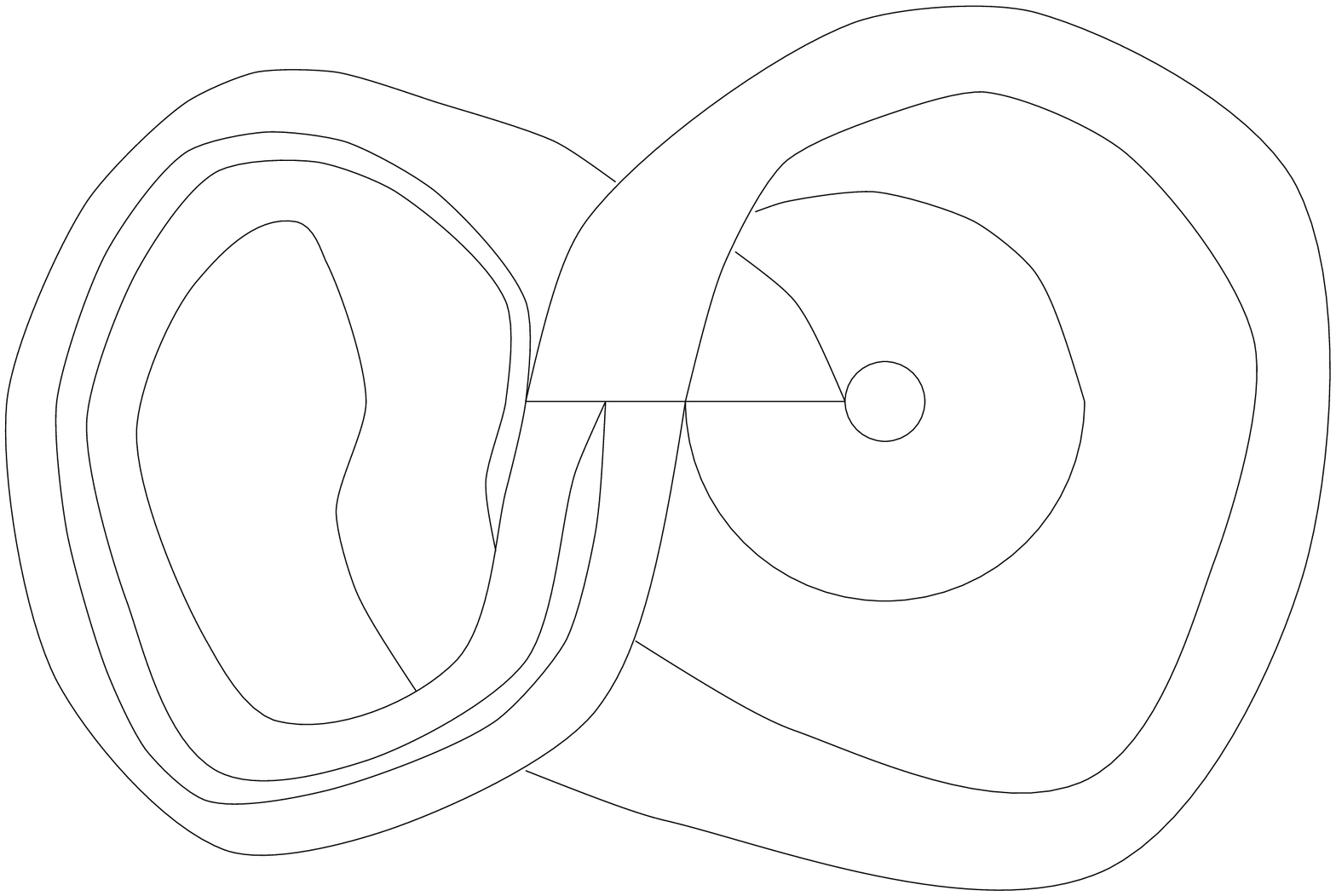} 

\pagebreak

k)      \hspace{2.5in}          l)

\includegraphics[height=1.5in,width=2in]{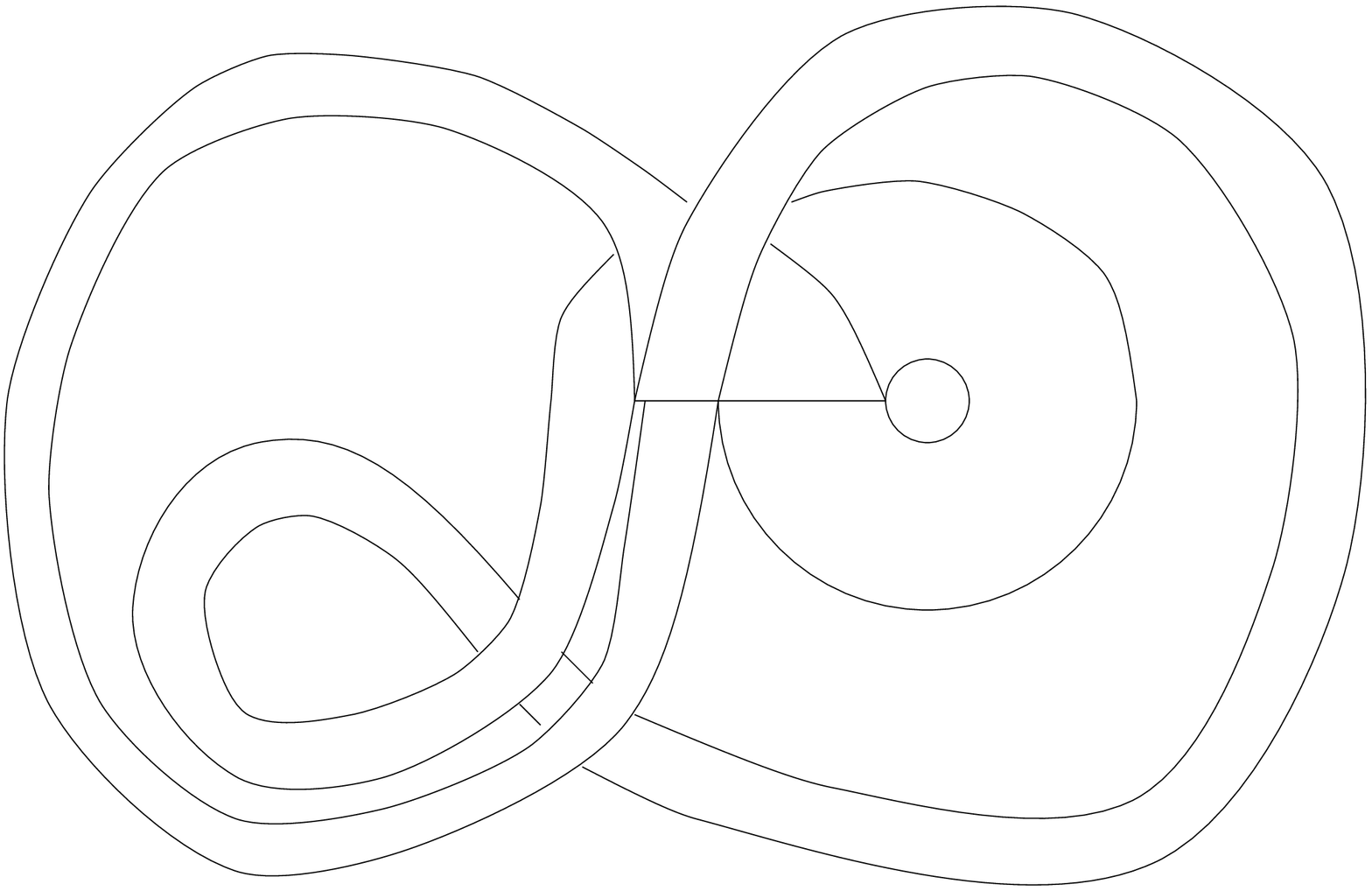} \hspace{.25in}
\includegraphics[height=1.5in,width=2in]{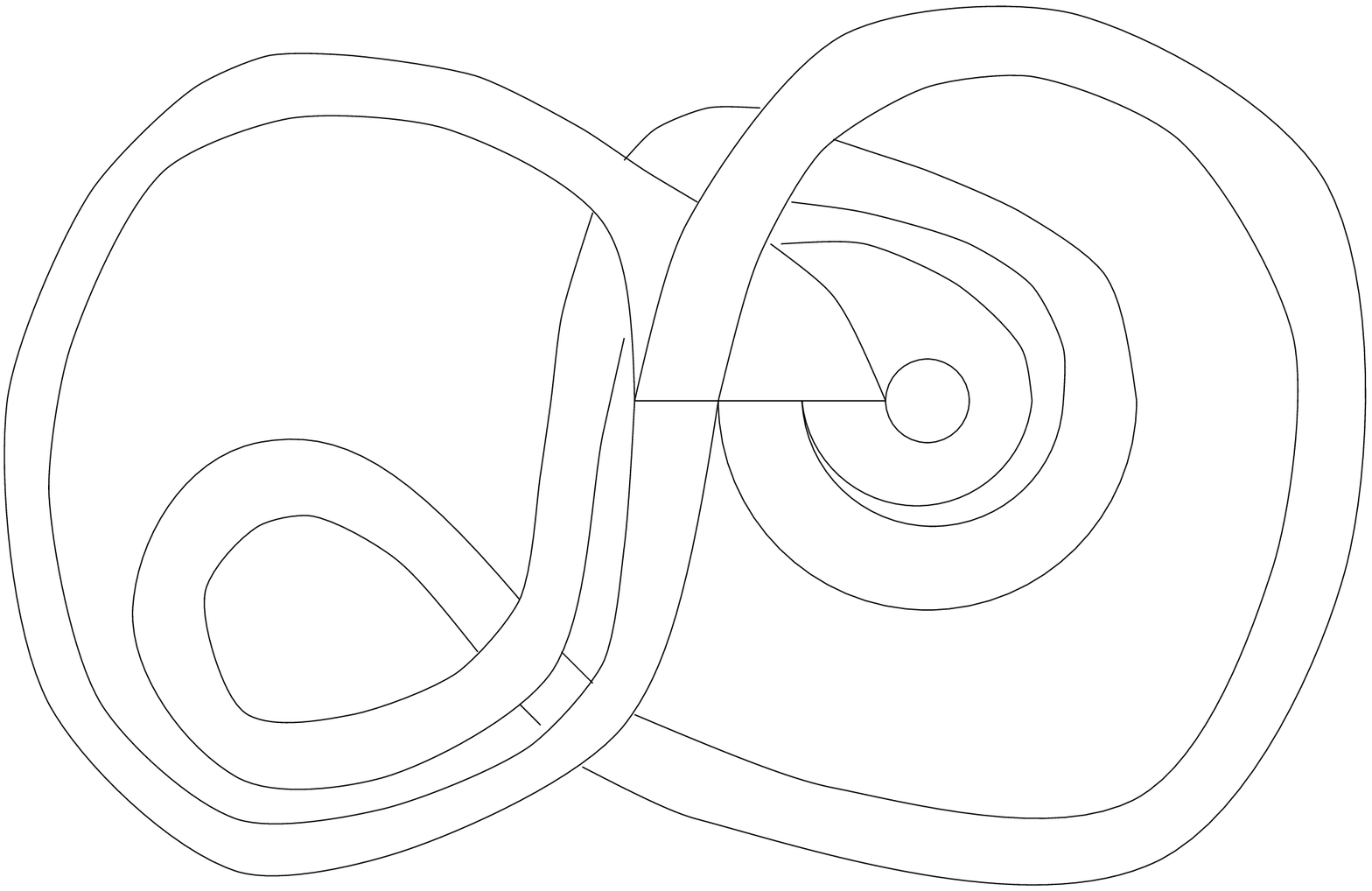} 

m)      \hspace{2.5in}          n)

\includegraphics[height=1.5in,width=2in]{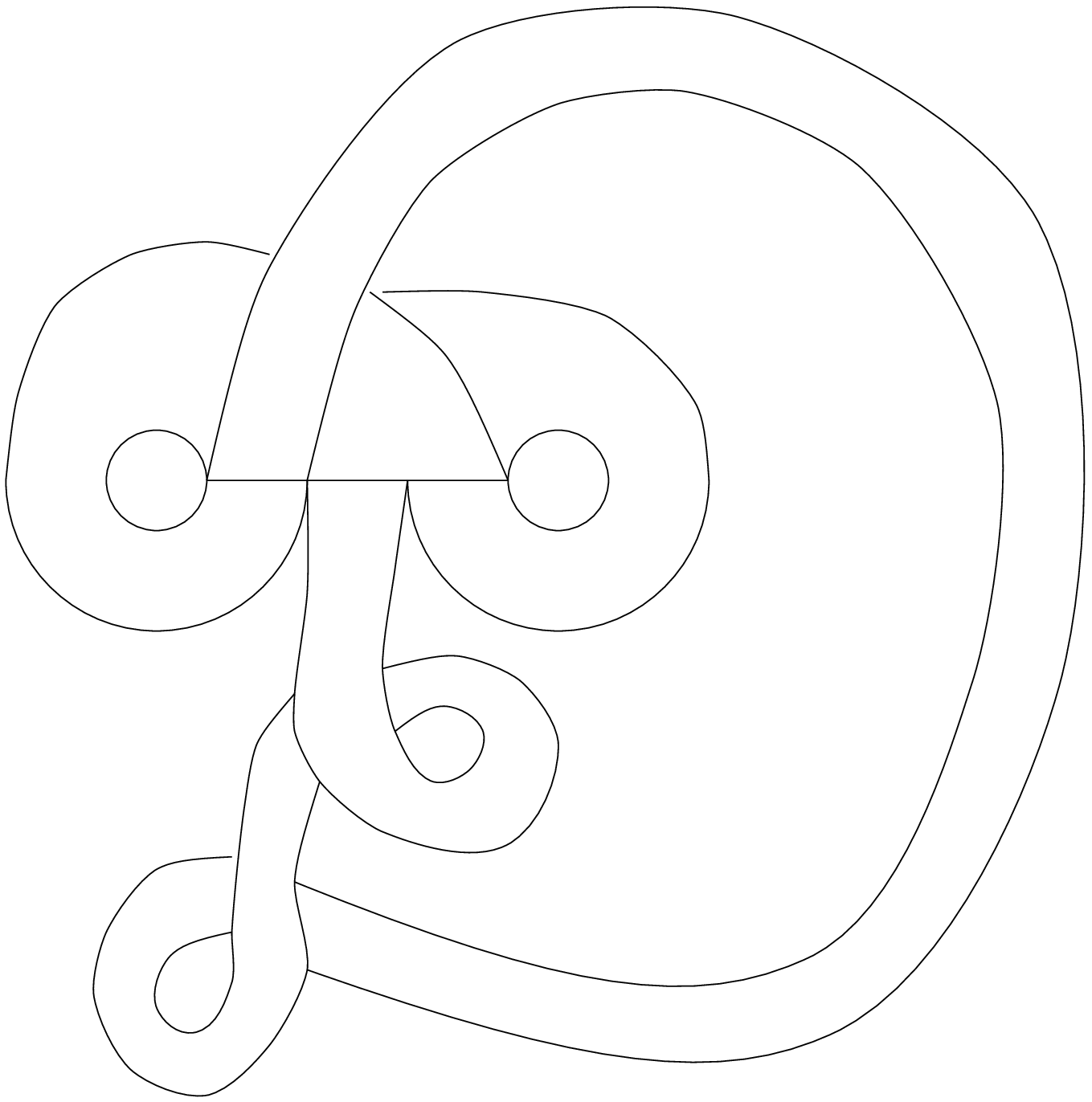} \hspace{.25in}
\includegraphics[height=1.5in,width=2in]{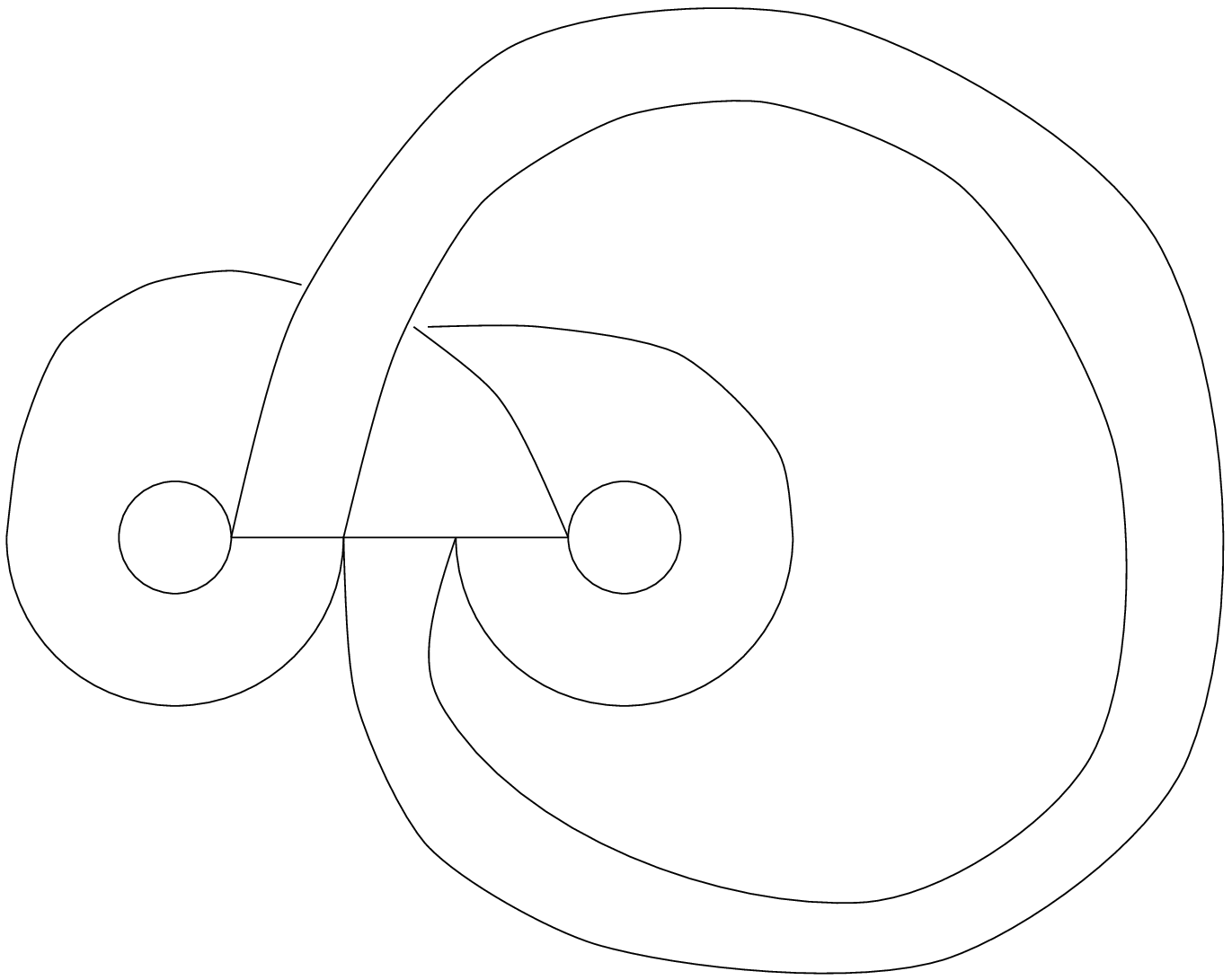}

\hspace*{1in} o)

\hspace*{1in}\includegraphics[height=1.5in,width=2in]{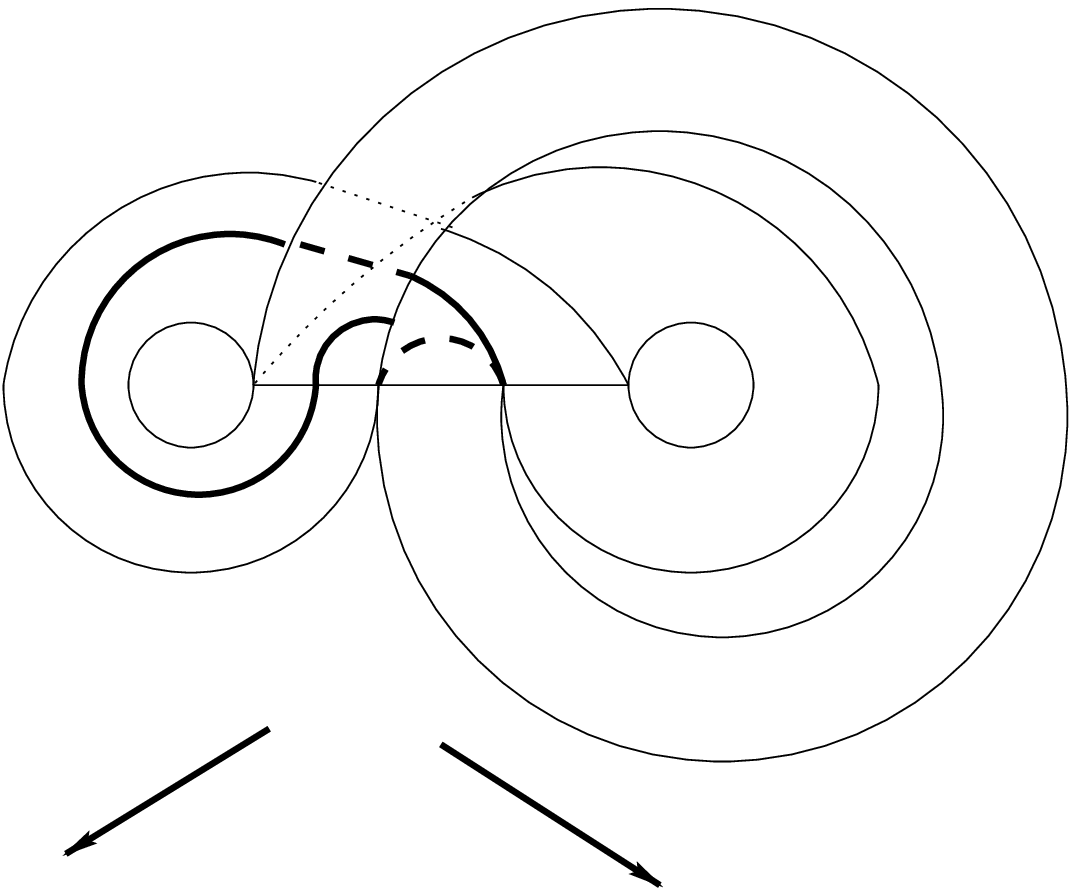} 

p)      \hspace{2.5in}          t)

\includegraphics[height=1.0in,width=1.5in]{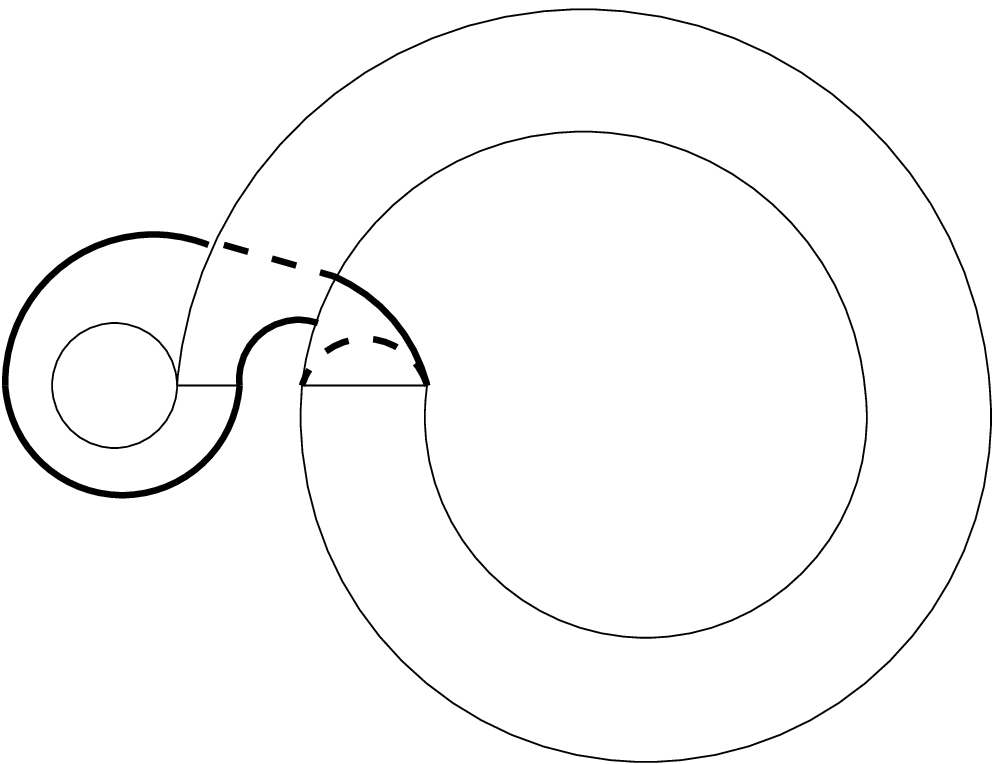} \hspace{1.25in}
\includegraphics[height=1.0in,width=1.5in]{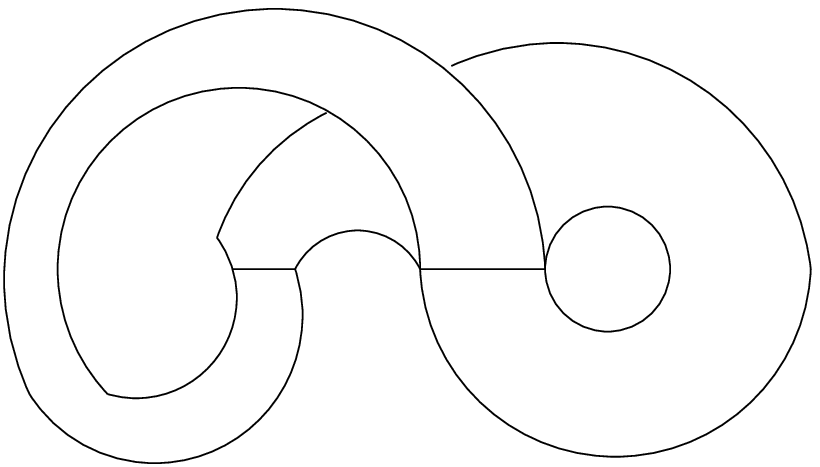}

q)      \hspace{2.5in}          u)

\includegraphics[height=1in,width=1.5in]{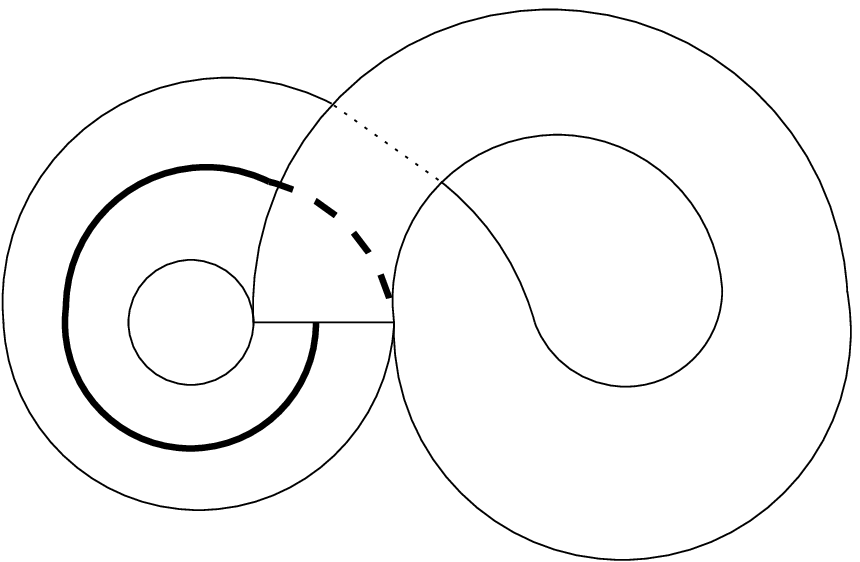} \hspace{1.25in}
\includegraphics[height=1in,width=1.5in]{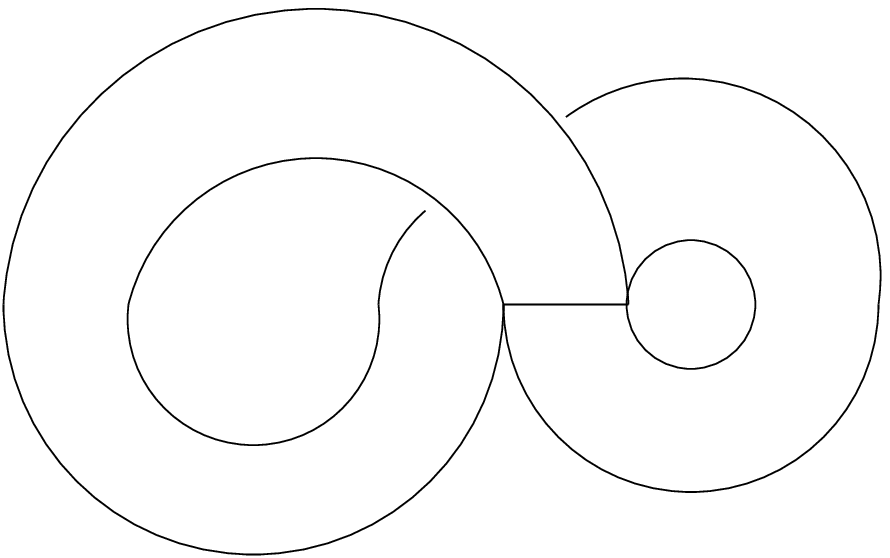} 

\pagebreak

r)      \hspace{2.5in}          v)

\includegraphics[height=1in,width=1.5in]{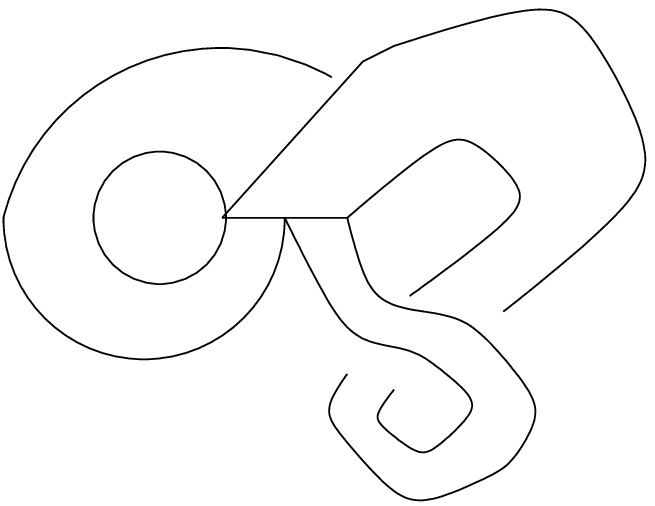} \hspace{1.25in}
\includegraphics[height=1in,width=1.5in]{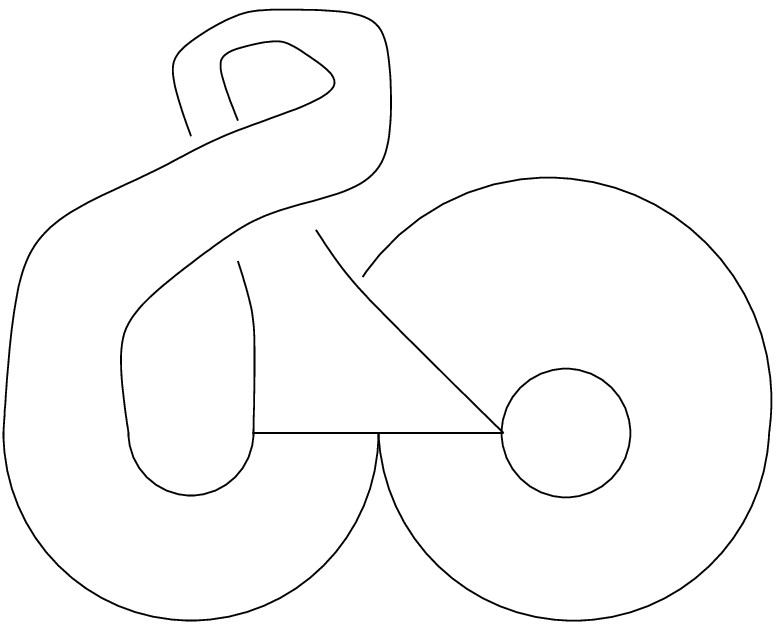} 

s)      \hspace{2.5in}          w)

\includegraphics[height=1in,width=1.5in]{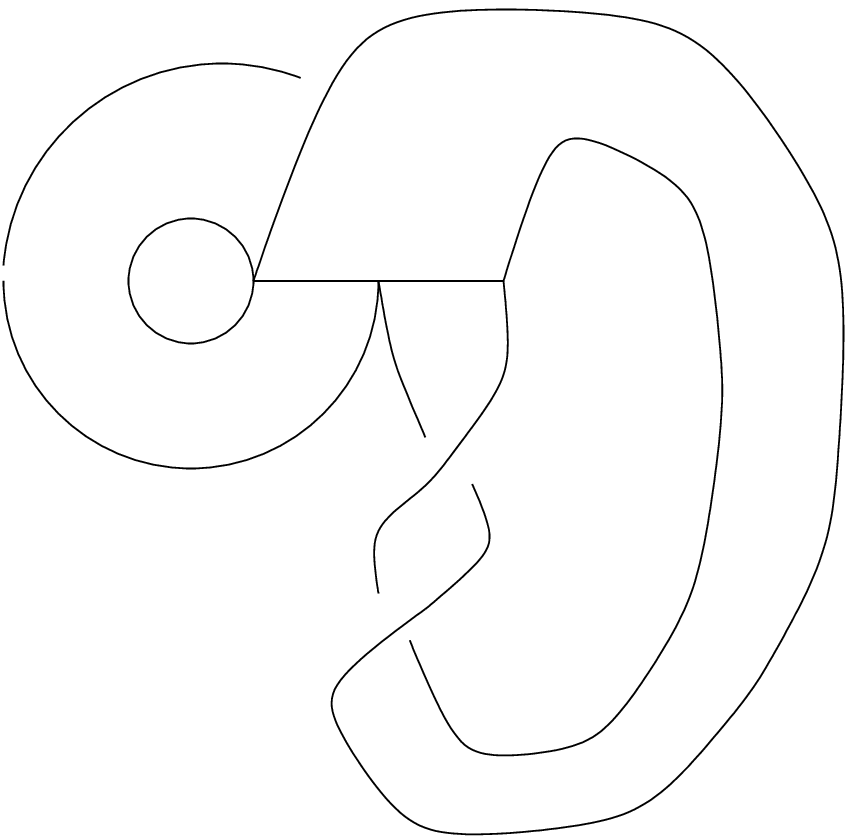} \hspace{1.25in}
\includegraphics[height=1in,width=1.5in]{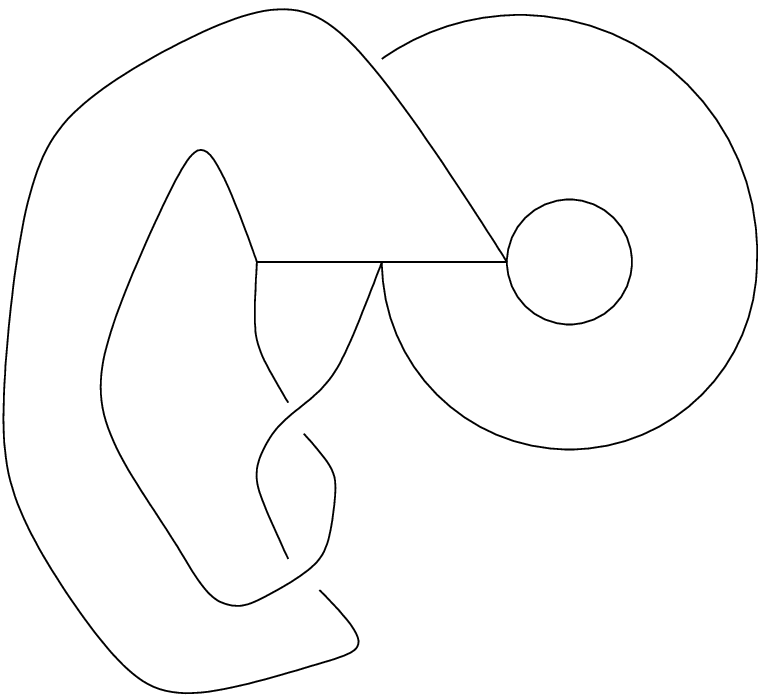} 

\begin{figure}[!h]
	\caption{Figures (a)-(w) prove Theorem~\ref{thm_main}}
	\label{a-w}
\end{figure}

% @@@@@@@@@@@@@@@@@@@@@@@ FIGURE @@@@@@@@@@@@@@@@@@@@@@@@@@@@@@@@@@@@@@@@@@
\begin{figure}[!htb]
	\begin{center} 
        \includegraphics[height=1.5in]{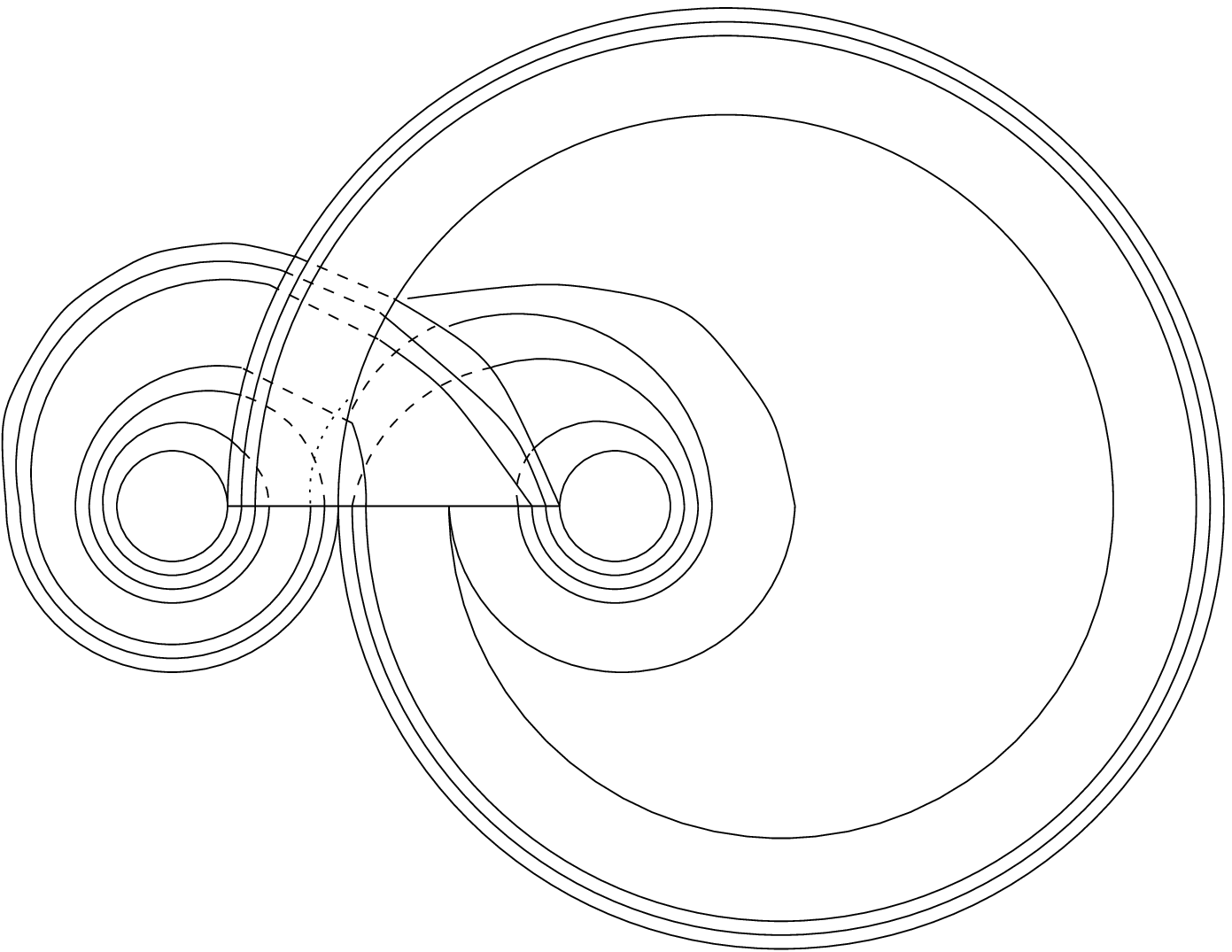}

	\includegraphics[height=2in]{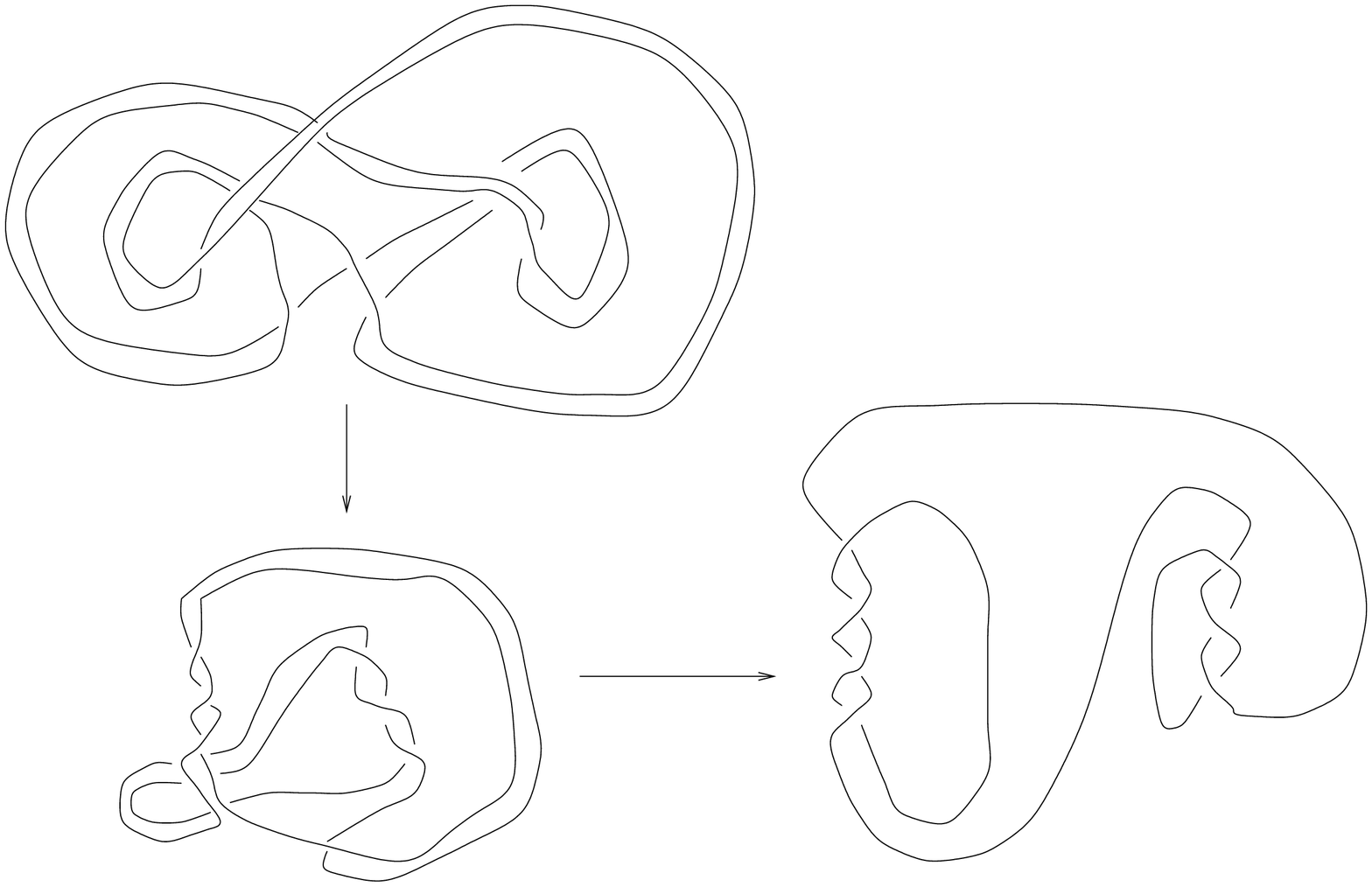}
	\end{center}
	\caption{Composite knot}
	\label{fig6.eps}
\end{figure}
% @@@@@@@@@@@@@@@@@@@@@@@@@@@@@@@@@@@@@@@@@@@@@@@@@@@@@@@@@@@@@@@@@@@@@@@@@@

% @@@@@@@@@@@@@@@@@@@@@@@ FIGURE @@@@@@@@@@@@@@@@@@@@@@@@@@@@@@@@@@@@@@@@@@
\begin{figure}[!htb]
	\begin{center} 
	\includegraphics[height=2in]{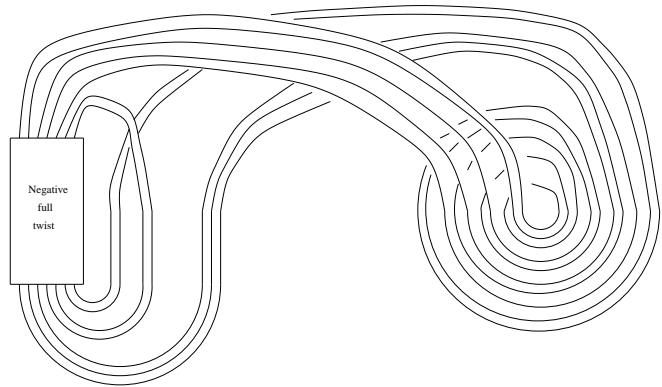}
	\end{center}
	\caption{Composite knot on $L(-2,0)$}
	\label{bigknot.eps}
\end{figure}
% @@@@@@@@@@@@@@@@@@@@@@@@@@@@@@@@@@@@@@@@@@@@@@@@@@@@@@@@@@@@@@@@@@@@@@@@@@

Thus {\sc Theorems} 1 and 2 imply that for all negative integers $n$,
$L(0,n)$ contains composite knots. In fact {\sc Theorem} 2 shows how to 
construct many examples of composite knots on $L(0,-2)$ and hence for any
Lorenz-like template with an even number of negative twists on one its 
branches. The theorem below does this for the odd cases.

\begin{thm}
	As sets of knots $L(0,-4)$ and $\tilde{L}(1,-(2n-1))$, for all $n>0$,
	are subsets of $L(0,-1)$.
\end{thm}

\begin{figure}[!h]
	\begin{center} 
	\includegraphics[height=3in]{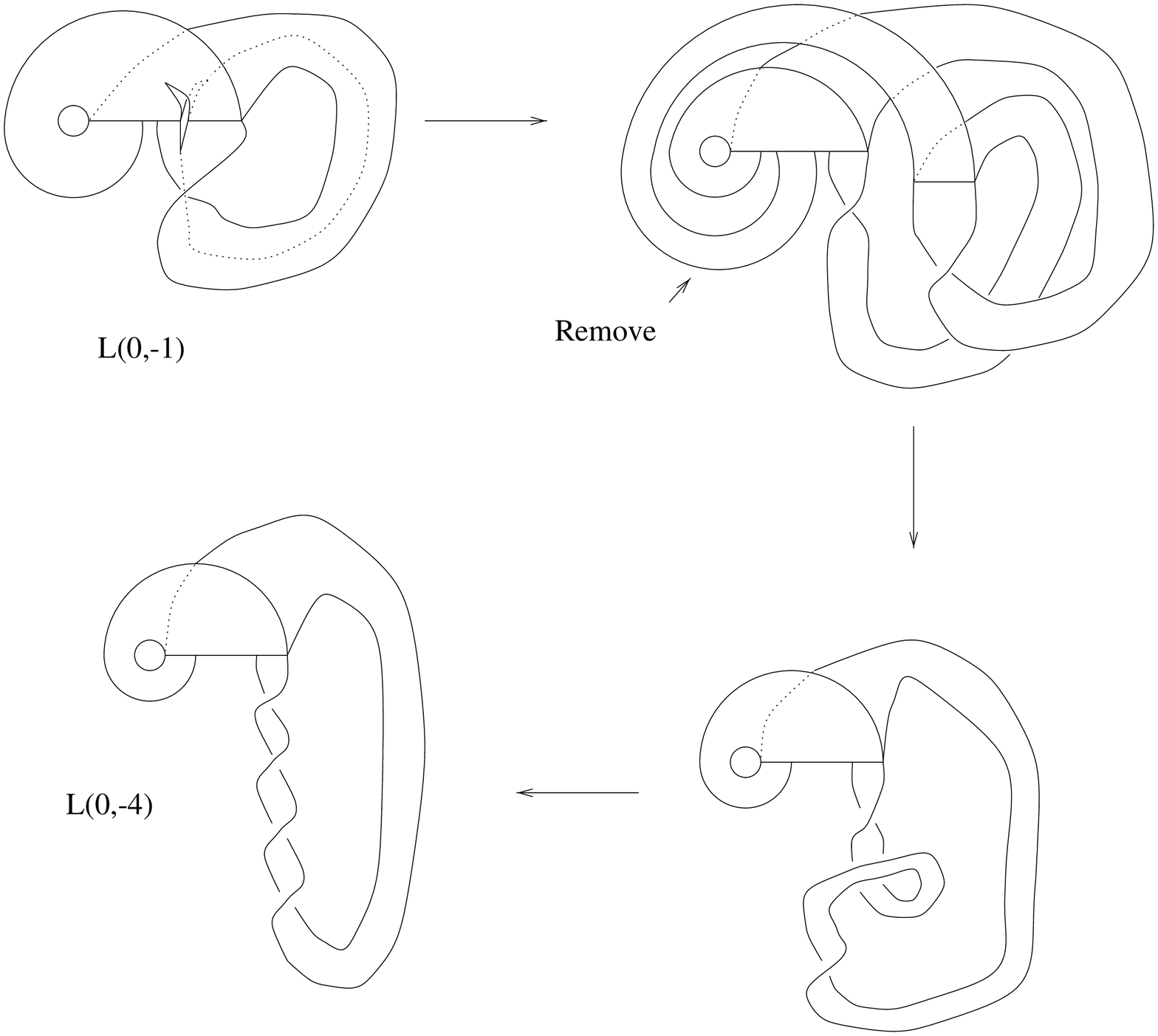}
	\end{center}
	\caption{$L(0,-4) \subset L(0,-1)$}
	\label{4in1.eps}
\end{figure}

\begin{figure}[htb]
	\begin{center} 
	\psfrag{L(0,-1)}{$\tilde{L}(0,-1)$}
	\psfrag{L(-1,1)}{$\tilde{L}(-1,1)$}
	\psfrag{L(-3,1)}{$\tilde{L}(-3,1)$}
	\psfrag{L(-5,1)}{$\tilde{L}(-5,1)$}
	\includegraphics[height=7in]{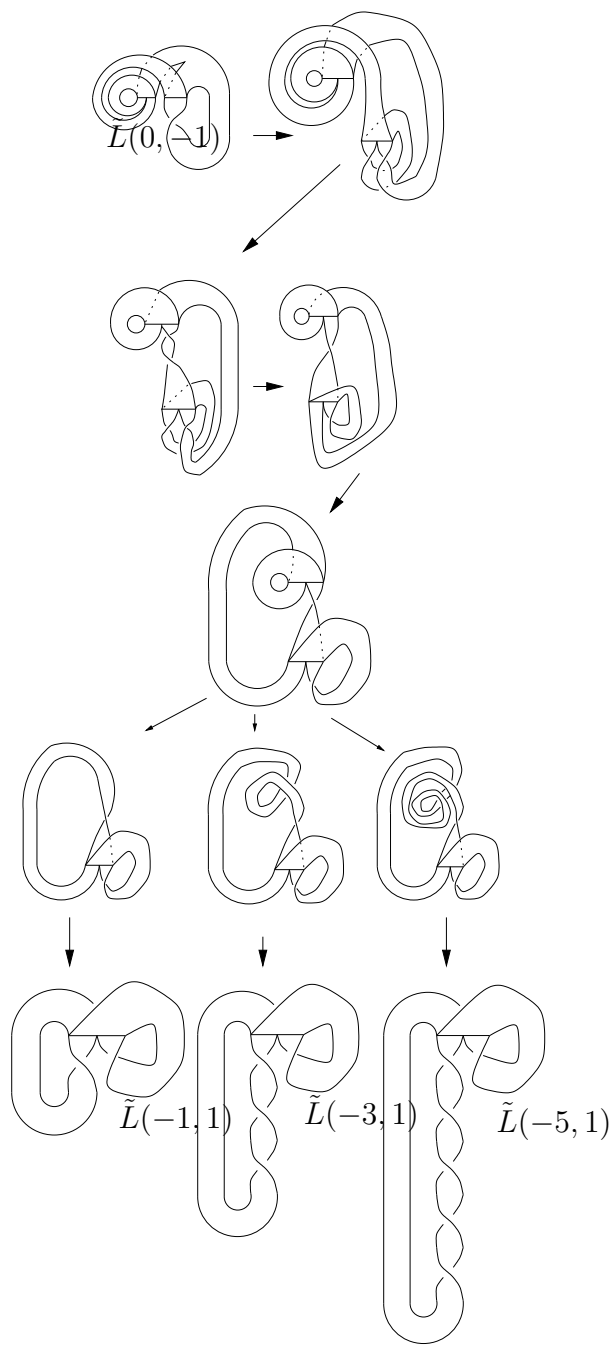}
	\end{center}
	\caption{Rest of Theorem 3}
	\label{rest.eps}
\end{figure}

\begin{proof}
Figure~\ref{4in1.eps} shows that the knots in $L(0,-4)$ are in $L(0,-1)$.
Figure~\ref{rest.eps} completes the proof.
\end{proof}

% @@@@@@@@@@@@@@@@@@@@@@@ FIGURE @@@@@@@@@@@@@@@@@@@@@@@@@@@@@@@@@@@@@@@@@@
% @@@@@@@@@@@@@@@@@@@@@@@@@@@@@@@@@@@@@@@@@@@@@@@@@@@@@@@@@@@@@@@@@@@@@@@@@@

% @@@@@@@@@@@@@@@@@@@@@@@ FIGURE @@@@@@@@@@@@@@@@@@@@@@@@@@@@@@@@@@@@@@@@@@
% @@@@@@@@@@@@@@@@@@@@@@@@@@@@@@@@@@@@@@@@@@@@@@@@@@@@@@@@@@@@@@@@@@@@@@@@@@

\clearpage

\noindent
{\bf Acknowledgement:} Supported in part by NSF grant no. DMS-9001973.

\end{document}